\documentclass[sn-mathphys]{sn-jnl}

\usepackage{anyfontsize}
\usepackage{amssymb}
\usepackage{cmap}
\usepackage{hyperxmp}
\usepackage{mathtools}
\usepackage{tikz-cd}
\usepackage{program}
\mathtoolsset{showonlyrefs}
\usepackage{enumerate}
\usepackage{todonotes}

\jyear{2025}%

\usepackage{booktabs,tabularx,siunitx,threeparttable}
\theoremstyle{thmstyleone}%

\theoremstyle{thmstyletwo}%
\newtheorem{rem}{Remark}%

\theoremstyle{thmstylethree}%

\newcommand{\RR}{\mathbb{R}} 

\newcommand{\norm}[1]{\left\Vert#1\right\Vert}          

\renewcommand{\im}{\operatorname{Im}} %

\begin{document}

\title[]{Weighted Hodge Laplacians on Manifolds with Boundary}

\author*[1]{\fnm{Zhe} \sur{Su}}\email{zhs0011@auburn.edu}
\author[2]{\fnm{Yiying} \sur{Tong}}
\author[3,4]{\fnm{Guo-Wei} \sur{Wei}} 
\affil[1]{\orgdiv{Department of Mathematics and Statistics}, \orgname{Auburn University}, \orgaddress{\state{AL 36849}, \country{USA}}}
\affil[2]{\orgdiv{Department of Computer Science and Engineering}, \orgname{Michigan State University}, \orgaddress{\state{MI 48824}, \country{USA}}}
\affil[3]{\orgdiv{Department of Mathematics}, \orgname{University of Georgia}, \orgaddress{\state{GA 30602}, \country{USA}}}
\affil[4]{\orgdiv{Department of Biochemistry and Molecular Biology}, \orgname{University of Georgia}, \orgaddress{\state{MI 30602}, \country{USA}}}


\abstract{The spectrum of the Hodge Laplacian on differential manifolds encodes rich topological and geometric information and thus provides a powerful tool for analyzing data on manifolds. However, the classical unweighted formulation is restricted in its ability to study data with varying local features. To address this limitation, we propose a weighted Hodge Laplacian framework for manifolds with boundary, both in theory and in computation, by incorporating a weight function on the manifold. Under appropriate boundary conditions, we formulate the corresponding weighted de Rham–Hodge theory, in which the kernel of the weighted Hodge Laplacian coincides with the weighted harmonic space, and remains isomorphic to the de Rham cohomology of the underlying manifold. The harmonic spectrum of the weighted Hodge Laplacian captures the global topological information, while its non-harmonic spectrum encodes the local geometric property induced by the weight. The proposed framework therefore enables the study of topological and geometric features of data on manifolds across varying weights, and in addition, allows local structure to be highlighted by choosing weights that emphasize regions of interest. We demonstrate the effectiveness of the proposed method through proof-of-principle experiments in protein flexibility analysis, and the results show its promise.}



\maketitle






\tableofcontents

\newpage

\section{Introduction}

The Hodge Laplacian on differential manifolds, which arises from the de Rham-Hodge theory, allows one to study the topological and geometric structures in a unified analytic framework through differential forms. Its spectrum carries rich information about the underlying manifold. When acting on differential $k$-forms, its kernel is equal to the space of harmonic $k$-forms (under certain boundary conditions with the presence of a boundary), which, through the Hodge isomorphism, also represents the de Rham cohomology. The kernel dimension is given by the $k$-Betti number, which counts independent $k$-dimensional holes. Its non-harmonic part, i.e., the sequence of non-zero eigenvalues, reflects geometric and algebraic features of the underlying manifold. These connections make the Hodge Laplacian useful in many fields. The spectrum of the Laplacian-Beltrami operator, i.e., the Hodge Laplacian on functions, in particular, has been widely used in spectral shape analysis for shape retrieval \cite{lian2011shape}, shape representation and matching \cite{reuter2006laplace, rustamov2007laplace}, mesh segmentation \cite{de2008hierarchical, reuter2010hierarchical}, and surface parameterization \cite{mullen2008spectral}. The operator has been used to study problems such as steady state, heat diffusion, and wave propagation.  

Parallel to the study of Hodge Laplacians, combinatorial Laplacians \cite{eckmann1944harmonische,horak2013interlacing, lim2020hodge} extend the graph Laplacian \cite{kirchhoff1847ueber} to simplicial complexes and also provide a connection between their kernels and the homologies of the underlying simplicial complex. The combinatorial Laplacians have been effective in areas such as signal processing \cite{roddenberry2022signal}, biomolecule analysis \cite{wang2020persistent,chen2022persistent,meng2021persistent,qiu2023persistent}, etc. Note that the term ``Hodge Laplacians" is sometimes used in this combinatorial context. However, this use of the term can be misleading. Although combinatorial Laplacians admit discrete ``Hodge decompositions", only Hodge Laplacians on differential manifolds lead to physically/geometrically relevant components, such as curl-free, divergence-free, and harmonic parts as in vector calculus. For manifolds with boundary, Hodge Laplacians admit meaningful boundary conditions, such as the normal (Dirichlet) and tangential (Neumann) boundary conditions, which are of significance in applications such as electromagnetism and fluid dynamics, where the boundary has physical or geometric implications. In addition, Hodge Laplacians encode the Riemannian metric and thus reflect the intrinsic geometry of the underlying manifold. These rich notions, however, do not emerge from relatively simple combinatorial Laplacians defined on point clouds.  For readers interested in the similarities and differences between Hodge Laplacians and combinatorial Laplacians, see \cite{ribando2024combinatorial}. 

Recently in topological data analysis (TDA), both combinatorial and Hodge Laplacians have been extended to the persistent setting to address limitations of persistent homology \cite{edelsbrunner2002topological, zomorodian2005computing, adler2010persistent}, a fundamental tool in TDA, which fails to capture the homotopic geometric shape changes in the filtration, leading to persistent combinatorial Laplacians \cite{wang2020persistent,lieutier2014persistent} and persistent Hodge Laplacians \cite{chen2021evolutionary,su2024persistent}, respectively. In combination with machine learning, particularly deep neural networks \cite{cang2017topologynet}, these persistent Laplacians have been shown to be powerful in characterizing complex data in biological problems \cite{wang2020persistent,zhao2020rham,meng2021persistent,wee2022persistent,qiu2023persistent,su2024persistent}. In particular, persistent combinatorial Laplacians have successfully predicted the emerging dominant SARS-CoV-2 variants BA4 and BA5 \cite{chen2022persistent}, which is one of the most impressive achievements of persistent Laplacians and TDA. The achievement has motivated the development of a variety of persistent Laplacians for other topological domains~\cite{wei2021persistent,wang2023persistent,liu2021persistent,chen2023persistent,jones2025persistent}. A survey of these persistent topological Laplacians is available \cite{wei2025persistent}. For TDA frameworks beyond persistent homology, we refer the reader to the review paper \cite{su2025topological}.

The proposed weighted Hodge Laplacians provide another refinement of the standard Hodge Laplacian. Motivated by the study of diffusion processes on manifolds, the drifting Laplacian for functions, also called the Bakry-Emery Laplacian \cite{bakry2006diffusions}, is obtained by incorporating a smooth weight function on the manifold via a modified codifferential, and has been studied in \cite{ma2010extension, xia2014inequalities, ma2008convexity, ma2009convex, lu2012eigenvalues} for its spectral properties. Its extension to differential forms is the drifting Hodge Laplacian. For non-compact manifolds, a de Rham-Hodge theory for drifting Hodge Laplacians is established in \cite{bueler1999heat}, where the results also apply to compact complete manifolds. The eigenvalue estimates of drifting Hodge Laplacians appeared in \cite{branding2024eigenvalue}. The study of drifting Hodge Laplacians has since become fundamental in the analysis of manifolds with density, diffusion process, and the study of Ricci flow. A related generalization of the standard Hodge Laplacian is the Witten Laplacian \cite{witten1982supersymmetry}, which was introduced in the context of supersymmetric quantum mechanics and Morse theory, and is closely connected to the topological quantum field theory \cite{atiyah1988topological} and recent developments in quantum topological data analysis (qTDA) \cite{di2024quantum}. In contrast to the drifting Hodge Laplacians, it modifies both the differential and codifferential operators. The Witten-Hodge theory for manifolds with boundary and the equivariant cohomology has been discussed in \cite{al2012witten}. These approaches enable the local analysis of data on manifolds, which is very useful in mathematical modeling of generative artificial intelligence (AI).   

Computationally, the discretization of Hodge Laplacians can be realized by the discrete exterior calculus (DEC) \cite{desbrun2006discrete} and the finite element exterior calculus (FEEC) \cite{arnold2018finite}. These methods are consistent with the continuous theory and thus preserve the correct cohomology of the underlying manifold, leading to accurate eigenvalue computation of the Hodge Laplacians. Based on these approaches, several computational frameworks have been developed for Hodge Laplacians and Hodge decompositions \cite{arnold2010finite,zhao20193d, ribando2024combinatorial,su2025topology}, where the underlying manifolds are represented either as simplicial meshes or as sublevel sets of level-set functions on regular Cartesian grids. In the latter case, a computational algorithm for persistent Hodge Laplacians is given in \cite{su2024persistent}. Despite the progress, these developments are unweighted, i.e., they do not incorporate weight functions.

The goal of this work is to develop a framework for weighted Hodge Laplacians on manifolds with boundary, both in theory and in computation. We work with the drifting Hodge Laplacian, which modifies only the codifferential. We formulate the corresponding weighted de Rham Hodge theory, in which, the kernels of the weighted Hodge Laplacians coincide with the (weighted) harmonic spaces under boundary conditions, and with the cohomology of the manifold.  The computational framework is implemented using DEC on a regular Cartesian grid with built-in boundary conditions. The weighted Hodge Laplacian can extract both the topological and geometric information from data across varying weights or densities, and in particular, allows one to highlight and differentiate \emph{local} features by choosing weight functions that emphasize regions of interest. To demonstrate the effectiveness of our framework, we build blind machine learning prediction models for the prediction of protein B-factors, which reflect the local flexibility of protein atoms, and compare them with the previous models. Experimental results show that our models achieve better performance, highlighting the promise of the proposed weighted Hodge Laplacian framework.

The rest of this article is organized as follows: In Sec.~\ref{sec:deRhamHodge}, we recall the de Rham-Hodge theory on manifolds with boundary, including the de Rham cohomology and the standard Hodge Laplacians. In Sec.~\ref{sec:deRhamHodge.weighted}, we explore the weighted setting by incorporating a smooth weight function on the manifold. We present the weighted de Rham cohomology and also introduce the weighted Hodge Laplacians. We also discuss how the kernels of Laplacians relate to the spaces of harmonic fields, and to the cohomology in this context. In Sec.~\ref{sec:discretization}, we provide a discretization framework for computing the weighted Hodge Laplacians on a regular Cartesian grid. Finally, we present experiments on protein flexibility analysis in Sec.~\ref{sec:proteinFlexibilityAnalysis}, and conclude the paper in Sec.~\ref{sec:conclusion}.

\section{De Rham Hodge Theory}
\label{sec:deRhamHodge}

In this section, we present the de Rham-Hodge theory on manifolds with boundary. The theory connects algebraic topology, differential geometry and analysis via differential forms, and it provides an intrinsic description of the global structure of differential manifolds through local differential operators. Under certain boundary conditions such as the normal (Dirichlet) or tangential (Neumann) boundary conditions, the kernel of Hodge Laplacian identifies with the space of harmonic forms and is naturally isomorphic to the relative or absolute cohomology of the underlying manifold. The Hodge Laplacian, therefore, is crucial in understanding the topological and geometric structure of a differential manifold.

We denote by $\Omega^k$ the space of all differential $k$-forms on $M$, where $M$ is a smooth, orientable and compact $m$-dimensional Riemannian manifold with boundary $\partial M$. The differential $d$, also called the exterior derivative, is a fundamental differential operator that maps $k$-forms to $(k\!+\!1)$-forms. It satisfies the Leibniz rule with respect to the wedge product $\wedge$ and the nilpotent property $dd = 0$. A form $\omega$ is called closed if $d\omega =0$, and exact if there is an $\eta\in\Omega^{k-1}$ such that $d\eta=\omega$.

A Riemannian metric $g$ on $M$ induces the Hodge star operator $\star$, which provides a canonical isomorphism between the space of $k$-forms and the space of $(m\!-\!k)$-forms. Using the Hodge star, one can then define the Hodge $L^2$-inner product as follows:
\begin{align}\label{eq:HodgeL2InnerProduct}
    (\omega, \eta) = \int_M\langle\omega, \eta\rangle_gd\mu = \int_M\omega\wedge\star\eta,
\end{align}
where $\langle\cdot, \cdot\rangle_g$ is the pointwise inner product metric on $\Omega^k$ induced by $g$ and $d\mu$ denotes the volume form induced by $g$.

The codifferential $\delta = (-1)^{m(k-1)+1}\star d\star $ maps $(k\!+\!1)$-forms to $k$-forms and also has the nilpotent property $\delta\delta = 0$. We then call a $k$-form $\omega$ coclosed if $\delta\omega=0$ and coexact if there is an $\eta\in\Omega^{k+1}$ such that $\omega = \delta\eta$. In the case of closed manifolds, i.e., manifolds without boundary, the codifferential $\delta$ can also be defined as the adjoint of the differential with respect to the Hodge $L^2$-inner product. However, in the presence of a boundary, this is no longer the case since the boundary term in general does not vanish, as shown in the formula below through integration by parts
\begin{align}
    (d\omega, \eta) = (\omega, \delta\eta) + \int_{\partial M}\omega\wedge\star\eta.
\end{align}
The adjointness of the codifferential $\delta$ and the differential $d$, however, is crucial, as it ensures the self-adjointness of the Hodge Laplacians, which guarantees real eigenvalues and the orthogonality of the eigenspaces, and additionally the connection between the kernels of Hodge Laplacians and the de Rham cohomology classes. These form the basis of the de Rham-Hodge theory.

\begin{rem}
The differential $d$ and the codifferential $\delta$ both generalize and unify the classical operators in vector calculus. For instance, in $\RR^3$, a $0$-form or a $3$-form can be identified with a scalar field, while a $1$-form or a $2$-form can be identified with a vector field. In this case, the differential $d$ applied to $0$-, $1$-, and $2$-forms yields the gradient operator $\nabla$, the curl operator $\nabla\times$, and the divergence operator $\nabla\cdot$, respectively. The codiffernetial $\delta$ corresponds to $-\nabla\cdot$, $\nabla\times$ and $-\nabla$ when applied to $1$-forms, $2$-forms, and $3$-forms, respectively. The nilpotent property $dd=0$ (equivalently $\delta\delta = 0$) leads to the vector calculus identities $\nabla\times\nabla = 0$ and $\nabla\cdot\nabla\times =0$. See \cite{desbrun2006discrete} for a detailed correspondence.
\end{rem}

To restore the adjointness of the differential $d$ and the codifferential $\delta$ in the presence of a boundary, one must impose certain boundary conditions. The most common choices are the normal (Dirichlet) and tangential (Neumann) boundary conditions. A differential form $\omega$ is normal if it vanishes on tangent vectors to the boundary, or tangential if the same condition holds for its Hodge dual $\star\omega$. Denote by $\Omega^k_n$ the set of normal differential $k$-forms and by $\Omega^k_t$ the set of tangential differential forms. We then have 
\begin{align} \label{eq.bc.n}
\Omega^k_n &= \{\omega\in\Omega^k\, \vert\, \omega\vert_{\partial M} = 0\};\\ \label{eq.bc.t}
\Omega^k_t &= \{\omega\in\Omega^k\, \vert\, \star\omega\vert_{\partial M} = 0\}.
\end{align}
The spaces $\Omega^k_n$ and $\Omega^{m-k}_t(M)$ are isomorphic under the Hodge star operator $\star$ by their definitions, which is known as the Hodge duality. In addition, the differential $d$ preserves the normal boundary conditions, while the codifferential $\delta$ preserves the tangential boundary conditions. This finally ensures that the de Rham cohomology is well-defined when restricting to the space of normal forms or the space of tangential forms. In the boundaryless case, both spaces $\Omega^k_n$ and $\Omega^k_t$ simplify to the entire space of differential forms $\Omega^k$.

\subsection{de Rham cohomology}

By restricting to the space of normal forms and the space of tangential forms, respectively, one can define the de Rham complex and its dual complex, which are sequences of the spaces of differential forms linked by $d$ and $\delta$ as follows:
\begin{align}\label{eq:deRhamcomplexes}
\cdots\overset{}{\longrightarrow} \Omega^{k-1}_n 
\overset{d^{k-1}}{\longrightarrow} &\Omega^{k}_n\overset{d^k}{\longrightarrow} \Omega^{k+1}_n \overset{}{\longrightarrow} \cdots \\
\cdots\overset{}{\longleftarrow} \Omega^{k-1}_t 
\overset{\delta^k}{\longleftarrow} &\Omega^{k}_t \overset{\delta^{k+1}}{\longleftarrow} \Omega^{k+1}_t \overset{}{\longleftarrow} \cdots.
\end{align}
Their cohomology modules then correspond to the relative de Rham cohomology $H^k_{dR}(M, \partial M)$ and the absolute de Rham cohomology $H^k_{dR}(M)$, i.e.,
\begin{align}\label{eq:cohomology}
H^k_{dR}(M, \partial M) &= \ker d^k / \im d^{k-1}\\
H^k_{dR}(M) &= \ker \delta^k / \im \delta^{k+1}.
\end{align}
Through Hodge duality and Lefschetz duality, the ranks of $H^k_{dR}(M, \partial M)$ and $H^k_{dR}(M)$ are given by the $(m\!-\!k)$-th Betti number $\beta_{m-k}$ and the $k$-th Betti number $\beta_k$ respectively. 

Note that here the domains of $d$ and $\delta$ in the definition of cohomology are restricted to the space of normal forms $\Omega^{k}_n(M)$ and the space of tangential forms $\Omega^{k}_t$, respectively. However, one can also equivalently define them by considering sequences of the entire space of forms $\Omega^k$ linked by $d$ or $\delta$. More details can be found in \cite{schwarz2006hodge}. These de Rham cohomology groups, by the de Rham theorem, are naturally isomorphic to the singular homology groups, and thus depend only on the topology of the underlying manifold and its boundary.

\subsection{Hodge Laplacians}
The $k$-th Hodge Laplacian $\Delta^k:\Omega^k\to\Omega^k$ is defined to be
\begin{align}\label{eq:HodgeLaplacian}
    \Delta^k = d^{k-1}\delta^k + \delta^{k+1}d^k.
\end{align}
For $0$-forms, $\Delta^0 = \delta^1d^0$ coincides with the usual Laplace-Beltrami operator for functions on the manifold. Its kernel $\ker\Delta^k$ is called the space of harmonic $k$-forms, which consists of all differential forms satisfying $\Delta^k\omega = 0$. 


The boundary constraint on the space of normal forms and tangential forms, however, is not enough to make the Hodge Laplacian self-adjoint and well-defined when its domain and codomain are restricted to these spaces. The issue can be addressed by imposing additional constraints. We thus restrict the Hodge Laplacian $\Delta$ to the following subspaces 
\begin{align} \label{eq.bc.n.subspace}
\bar\Omega^k_n &= \{\omega\in\Omega^k\, \vert\, \omega\vert_{\partial M} = 0,\, \delta\omega\vert_{\partial M} = 0\};\\ \label{eq.bc.t.subspace}
\bar\Omega^k_t &= \{\omega\in\Omega^k\, \vert\, \star\omega\vert_{\partial M} = 0,\, \star d\omega\vert_{\partial M} = 0\}.
\end{align}
Through Hodge duality,  these subspaces are pairwise isomorphic, i.e., $\bar\Omega^k_n\overset{\star}{\cong}\bar\Omega^{m-k}_t$. In addition, the following identity can be obtained after imposing these boundary conditions
\begin{align}\label{eq:LaplacianIdentity}
    (\Delta\omega, \omega) = (d\omega, d\omega) + (\delta\omega, \delta\omega).
\end{align}
Let $\mathcal{H}^k = \ker d^k\cap\ker\delta^k$ be the space of differential forms that are both closed and coclosed, also known as the space of harmonic fields. Denote by $\Delta^k_n$ and $\Delta^k_t$ the restrictions of the Hodge Laplacian $\Delta^k$ to $\bar\Omega^k_n$ and $\bar\Omega^k_t$, respectively. The following identifications are then a direct result of \eqref{eq:LaplacianIdentity}
\begin{align}
    \ker\Delta^k_n = \mathcal{H}^k_n,\quad \ker\Delta^k_t = \mathcal{H}^k_t,
\end{align}
where $\mathcal{H}^k_n = \mathcal{H}^k\cap\bar\Omega^k_n$ and $\mathcal{H}^k_t = \mathcal{H}^k\cap\bar\Omega^k_t$ are restrictions of the space of harmonic fields to normal and tangential differential forms, respectively. In addition, these two subspaces are isomorphic to the relative and the absolute de Rham cohomology groups \cite{friedrichs1955differential}, i.e., $\mathcal{H}^k_n\cong H^k_{dR}(M,\partial M)$ and $\mathcal{H}^k_t\cong H^k_{dR}(M)$. Under Lefschetz duality $H^k_{dR}(M)\cong H_k(M,\partial M)$ and Hodge duality $H^k_n\cong \star H^k_n=H^{m-k}_t$, $\dim H^k_n$ and $\dim H^k_t$ are therefore finite and given by the Betti numbers $\beta_{m-k}$ and $\beta_k$, respectively. Note that in general, the space of harmonic fields is only a subset of the space of harmonic forms in the presence of a boundary, i.e., $\mathcal{H}^k\subset\ker\Delta^k$, and is infinite dimensional \cite{schwarz2006hodge}. Thus, it does not provide accurate information about the manifold topology. The situation becomes simpler when the boundary vanishes. In this case, the subspaces $\mathcal{H}^k_n$ and $\mathcal{H}^k_t$ both coincide with $\mathcal{H}^k$, which is then equal to the space of harmonic forms $\ker\Delta$ and is isomorphic to the de Rham cohomology group $H^k_{dR}(M)$.

\begin{rem}

A central component in de Rham Hodge theory is the Hodge decomposition. In the case of manifolds with boundary, there is a Hodge-Morrey decomposition~\cite{morrey1956variational} given as:
\begin{align}\label{eq.morreyDecomp}
\Omega^k = d\Omega^{k-1}_n\oplus\delta\Omega^{k+1}_t(M) \oplus \mathcal{H}^k.
\end{align}
This $L^2$-orthogonal decomposition reduces to the classical Hodge decomposition in the boundaryless case that $M$ is a closed manifold. Its counterpart in vector calculus is often referred to as the Helmholtz-Hodge decomposition on compact domains in $\RR^2$ or $\RR^3$. In addition, there are $4$-component Hodge-Morrey-Friedrichs orthogonal decompositions~\cite{friedrichs1955differential}, and in particular, for compact domains in Euclidean spaces, one obtains a $5$-component Hodge decomposition \cite{shonkwiler2009poincare} as follows
\begin{align}\label{eq.hd.5subspaces0}
	\Omega^k = d\Omega^{k-1}_n\oplus\delta\Omega^{k+1}_t\oplus \mathcal{H}^{k}_n \oplus \mathcal{H}^{k}_t\oplus (d\Omega^{k-1}\cap\delta\Omega^{k+1}).
\end{align}
Due to the correspondence between one-forms and vector fields in Euclidean spaces, this $5$-component Hodge decomposition has been implemented for analyzing vector fields on surface triangle meshes, tetrahedral meshes \cite{poelke2016boundary,poelke2017hodge,zhao20193d,razafindrazaka2019consistent}, and Cartesian grids \cite{su2024hodge, su2025topology}. For an elementary exposition of this $5$-component Hodge decomposition in $\RR^3$ in terms of vector and scalar fields, see \cite{cantarella2002vector}.

\end{rem}

\section{Weighted de Rham-Hodge Theory}
\label{sec:deRhamHodge.weighted}

We now present the weighted de Rham-Hodge theory for manifolds with boundary, which generalizes the standard de Rham-Hodge theory by incorporating a smooth function $f\in C^{\infty}(M)$ defined on the manifold. There are several ways to define weighted Hodge Laplacians. Here we adopt the formulation of drifting Hodge Laplacians \cite{branding2024eigenvalue}, obtained by modifying only the codifferential operator. Specifically, we set the weighted codifferential $\delta_f:= e^f\delta e^{-f}$, which still maps $(k\!+\!1)$-forms to $k$-forms. It is straightforward to verify that $\delta_f$ is nilpotent, i.e., $\delta_f\delta_f = 0$, and it preserves the tangential boundary conditions. A form $\omega$ is then called $f$-coclosed if $\delta_f\omega = 0$, or $f$-coexact if there is $\eta\in\Omega^{k+1}(M)$ such that $\omega = \delta_f\eta$. 

A corresponding weighted measure on $M$ can be defined as $d\mu_f = e^{-f}d\mu$, where $d\mu$ is the volume form induced by $g$. The manifold $(M, g, d\mu_f = e^{-f}d\mu)$ in this case is often referred to as a Bakry-\'Emery manifold \cite{bakry2006diffusions}. One can then define a generalized $L^2$-inner product on the space of differential forms $\Omega^k$ by introducing the weighted measure, given as follows
\begin{align}\label{eq:WeightedL2InnerProduct}
    (\omega, \eta)_f = \int_M\langle\omega, \eta\rangle_ge^{-f}d\mu.
\end{align}
It coincides with a scalar multiple of the standard Hodge $L^2$-inner product~\eqref{eq:HodgeL2InnerProduct} when $f$ is constant. 

In the case of closed manifolds, one can easily verify that the weighted codifferential $\delta_f$ is adjoint to the differential $d$ with respect to the generalized $L^2$-inner product~\eqref{eq:WeightedL2InnerProduct}. The adjointness of $\delta_f$ and $d$ can be generalized to the case of manifolds with boundary by imposing the normal and tangential boundary conditions. To be specific, $\delta_f$ is adjoint to $d$ with respect to \eqref{eq:WeightedL2InnerProduct} when restricted to the space of ($f$-modified) normal forms and the space of tangential forms in the presence of a boundary, as detailed below.

\subsection{Weighted de Rham cohomology}

Note that in our weighted setting, the differential $d$ remains unchanged, only the codifferential is modified to incorporate the (log) weight $f$. Therefore, the de Rham complex and the associated relative de Rham cohomology are the same as in the unweighted case. Without loss of generality, we focus only on the dual complex defined by the weighted codifferential $\delta_f$ and the corresponding absolute de Rham cohomology. Since multiplying a differential form by a scalar function does not alter the (homogeneous) boundary conditions, the weighted codifferential $\delta_f$ preserves the tangential boundary conditions. The following dual complex is thus well-defined
\begin{align}\label{eq:deRhamcomplexes.weighted}
\cdots\overset{}{\longleftarrow} \Omega^{k-1}_t 
\overset{\delta^k_f}{\longleftarrow} &\Omega^{k}_t \overset{\delta^{k+1}_f}{\longleftarrow} \Omega^{k+1}_t \overset{}{\longleftarrow} \cdots.
\end{align}
We can then define a weighted version of the absolute de Rham cohomology  
\begin{align}\label{eq:cohomology.weighted}
H^k_{dR}(M, f) &= \ker \delta^k_f / \im \delta^{k+1}_f,
\end{align}
where the domain of $\delta^k_f$ is restricted to the space of tangential forms. Although the weight $f$ appears in the definition of the weighted absolute de Rham cohomology~\eqref{eq:cohomology.weighted}, this cohomology is isomorphic to the standard absolute de Rham cohomology in the unweighted case \eqref{eq:cohomology}, i.e.,
\begin{align}
    H^k_{dR}(M, f)\cong H^k_{dR}(M),
\end{align}
and thus is independent of the choice of $f$, and its dimension is given by the $k$-th Betti number $\beta_k$. This is an immediate result of the chain map $T:(\Omega^*_t, \delta_f)\to(\Omega^*_t, \delta), T(\omega) = e^{-f}\omega$ between the two complexes. The cohomology group~\eqref{eq:cohomology.weighted} is therefore fully determined by the topology of the underlying manifold.

\subsection{Weighted Hodge Laplacians}

Using the weighted codifferential $\delta_f$, one can define a weighted Hodge Laplacian (drifting Hodge Laplacian) $\Delta _f:\Omega^k\to\Omega^k$ as follows
\begin{align}\label{eq:weightedHodgeLaplacian}
\Delta^k_f:=d^{k-1}\delta^k_f+\delta^{k+1}_f d^k. 
\end{align}
When $f$ is constant, the weighted Hodge Laplacian~\eqref{eq:weightedHodgeLaplacian} coincides with the standard Hodge Laplacian~\eqref{eq:HodgeLaplacian}. The kernel of $\Delta^k_f$, denoted by $\ker{\Delta^k_f}$, is the space of weighted harmonic forms on $M$. 

The de Rham-Hodge theory for the drifting Hodge Laplacian \eqref{eq:weightedHodgeLaplacian} on non-compact manifolds has been established in \cite{bueler1999heat}. Their results also apply to compact complete manifolds. A characterization of the eigenvalues of the drifting Hodge Laplacian on manifolds immersed in $\RR^n$ can be found in \cite{branding2024eigenvalue}.

\begin{rem}
    In the case of degree $k = 0$, the weighted Hodge Laplacian $\Delta^k_f$ is known as the Bakry-\'Emery Laplacian \cite{bakry2006diffusions} given by
\begin{align}
    \Delta^0_f\phi= \Delta^0\phi - g(\nabla_f, \nabla\phi),
\end{align}
where $\phi$ is a function defined on $M$. It reduces to the Laplace-Beltrami operator when $f$ is constant.
\end{rem}

\begin{rem}
    In addition to the drifting Hodge Laplacian \eqref{eq:weightedHodgeLaplacian}, another weighted Hodge Laplacian, called Witten Laplacian \cite{witten1982supersymmetry}, can be obtained by modifying both the differential and codifferential given by
    \begin{align}
        \bar\Delta^k_f := d^{k-1}_f\delta^k_f+\delta^{k+1}_f d^k_f,
    \end{align}
    which was first introduced in the study of supersymmetry and Morse theory. Here the weighted differential and codifferential are defined by $d_f = e^{-f}de^{f}$ and $\delta_f = e^f\delta e^{-f}$, respectively. Note that Witten Laplacian remains in the same analytic framework as the standard Hodge Laplacian. All standard results, such as the adjointness of $d_f$ and $\delta_f$, the self-adjointness of Laplacian, and Hodge decomposition, are with respect to the standard Hodge inner product \eqref{eq:HodgeL2InnerProduct}. There is no need to redefine the $L^2$-inner product \eqref{eq:WeightedL2InnerProduct} in this context. The Witten-Hodge theory for manifolds with boundary and equivariant cohomology has been studied in \cite{al2012witten}.
\end{rem}

Before establishing the identifications between the kernels of the weighted Hodge Laplacians, the corresponding harmonic spaces, and the cohomology groups, we introduce the following subspace of normal forms adapted to $\delta_f$:
\begin{align}\label{eq.bc.f.n.subspace}
\bar\Omega^k_{n,f} &= \{\omega\in\Omega^k\, \vert\, \omega\vert_{\partial M} = 0,\, \delta_f\omega\vert_{\partial M} = 0\}.
\end{align}
The map $\star_f = e^{-f}\star$, called the weighted Hodge star, provides an isomorphism between this space of weighted normal forms $\bar\Omega^k_{n,f}$ and the space of tangential forms $\bar\Omega^{m-k}_t$, i.e., $\bar\Omega^k_{n,f}\overset{\star_f}{\cong}\bar\Omega^{m-k}_t$.
The claim follows from the fact that the Hodge star $\star$ is an isomorphism between $k$-forms and $(m\!-\!k)$-forms, and the functional multiple $e^{-f}>0$ never vanishes, and, in addition, the weighted Hodge star provides a bijection between the two subspaces.
Through integration by parts, the following identity can then be obtained when restricting $\Delta^k_f$ to the subspaces $\bar\Omega^k_{n,f}$ and $\bar\Omega^k_t$:
\begin{align}\label{eq:LaplacianIdentity.weighted}
(\Delta_f\omega, \omega)_f = \left(d\omega, d\omega \right)_f + \left(\delta_f\omega, \delta_f\omega \right)_f.
\end{align}
Note that $\delta_f =  e^f\delta e^{-f}=e^f [(-1)^k\star^{-1} d \star] e^{-f}=(-1)^k\star_f^{-1}d\star_f.$

Let $\mathcal{H}^k_f = \ker d\cap\ker\delta_f$ be the space of weighted harmonic fields, consisting of the differential forms that are both closed and $f$-coclosed. This space is infinite-dimensional and is only a subspace of weighted harmonic forms $\ker{\Delta_f}$, analogous to the unweighted case. We thus focus on its subspaces $\mathcal{H}^k_{n,f} = \mathcal{H}^k_f\cap\bar\Omega^k_{n,f}$ and $\mathcal{H}^k_t = \mathcal{H}^k_f\cap\bar\Omega^k_t$. Denote by $\Delta^k_{n,f}$ and $\Delta^k_{t,f}$ the restrictions of the weighted Hodge Laplacian $\Delta^k_f$ to $\bar\Omega^k_{n,f}$ and $\bar\Omega^k_t$. It follows directly from the identity~\eqref{eq:LaplacianIdentity.weighted} that
\begin{align}
    \ker\Delta^k_{n,f} = \mathcal{H}^k_{n,f},\quad \ker\Delta^k_{t,f} = \mathcal{H}^k_{t,f}.
\end{align}
In addition, there is an isomorphism given by the weighted Hodge star:
\begin{align}
    \mathcal{H}^k_{n,f}\overset{\star_f}{\cong}\mathcal{H}^{m-k}_{t, -f}.
\end{align}
These two subspaces are also isomorphic to the relative and the absolute de Rham cohomology groups, i.e., $\mathcal{H}^k_{n,f}\cong H^k_{dR}(M,\partial M)$ and $\mathcal{H}^k_{t, f}\cong H^k_{dR}(M)$, which follows by adapting the classical Hodge-theoretical arguments to the weighted setting under the respective boundary conditions. Therefore, they are finite-dimensional with dimensions given by the Betti numbers $\beta_{m-k}$ and $\beta_k$, respectively. The spaces $\mathcal{H}^k_{n,f}$ and $\mathcal{H}^k_{t,f}$ are fully determined by the manifold topology.

These identifications allow us to study the manifold topology through the weighted Hodge Laplacians on weighted normal and tangential differential forms. Note that while the standard Hodge Laplacians encode both the global topology and the local geometry of the manifold determined by the metric, introducing the weight $f$, in addition, allows us to study features that vary across weights and, in particular, highlight and differentiate the local geometric features in a uniform setting by selecting weights that emphasize regions of interest.

\begin{rem}
    A direct consequence of the adjointness of $d$ and $\delta_f$ with respect to the $L^2$-inner product~\eqref{eq:WeightedL2InnerProduct} establishes a weighted Hodge decomposition. It thus provides a data-driven approach to encode the local geometric information in the decomposition. In this context, the space of differential forms can be orthogonally decomposed with respect to \eqref{eq:WeightedL2InnerProduct} as follows 
\begin{align}
    \Omega^k = d\bar\Omega^{k-1}_{n,f}\oplus_f\delta_f\bar\Omega^{k+1}_t\oplus_f\mathcal{H}^k_f.
\end{align}
In addition, one naturally expects a weighted $4$-component Hodge-Morrey-Friedrichs decompositions with respect to the inner product~\eqref{eq:WeightedL2InnerProduct} in analogy with the unweighted case, and a weighted $5$-component Hodge decomposition when $M$ is a compact domain in a Euclidean space:
\begin{align}\label{eq.hd.5subspaces}
	\Omega^k = d\bar\Omega^{k-1}_{n,f}\oplus_f\delta_f\bar\Omega^{k+1}_t\oplus_f \mathcal{H}^{k}_{n,f} \oplus_f\mathcal{H}^{k}_{t, f}\oplus (d\Omega^{k-1}\cap\delta_f\Omega^{k+1}).
\end{align}
A proof of these decompositions is beyond the scope of this paper, and we will not address it further.
\end{rem}

\section{Discretization of the Laplacians}
\label{sec:discretization}

We now present a discretization of the weighted Hodge Laplacian for compact domains in low-dimensional Euclidean spaces, and also propose a weighted version of the boundary-induced graph (BIG) Laplacian \cite{ribando2024combinatorial} to facilitate computations. This framework builds on the discretization developed for the standard Hodge Laplacians \cite{su2024persistent} for regular Cartesian grids, where discrete exterior calculus (DEC) is used to discretize all differential operators and differential forms. The DEC-based approach enables efficient and accurate numerical algorithms that rely solely on matrix algebra and are consistent with the smooth case, preserving the topology of the manifold. 

As in \cite{su2024persistent}, we model the underlying manifold as a sublevel set of a level set function defined on a regular Cartesian grid, which is also referred to as the Eulerian formulation of a manifold. In this setting, all vertices, edges, faces, and cells are fixed, which greatly simplifies the data structure and discrete differential operators compared to the Lagrangian formulation, where the manifold is discretized by a triangular or a tetrahedral mesh. The Eulerian formulation enables the construction of the BIG Laplacians, which avoid the use of Hodge star operators while still capturing the correct cohomology. This simplification makes the computation highly efficient.

\subsection{Discretization on the entire grid}

Denote by $I_m$ an $m$-dimensional regular Cartesian grid with grid spacing $l$, whose cells are oriented consistently with the coordinate axes. We refer to this as the primal grid. Its dual, called the dual grid, is the staggered grid with grid points being the centers of the primal $m$-cells of $I_m$. For each primal $k$-cell $\sigma_k$, the dual $(m\!-\!k)$-cell is formed by the dual grid points associated with the primal $m$-cells incident to $\sigma_k$.

A discretization of a differential $k$-form $\omega$ on $I_m$, following the de Rham map, is simply a list of values with each given by its integral over the corresponding oriented primal $k$-cell. The discrete differential $D^I_k$ on $I_m$ encodes the signed incidence between the primal $(k\!+\!1)$-cells and $k$-cells. Specifically, it is the transpose of the cell boundary operator $\partial_{k+1}$ on $(k\!+\!1)$-cells following Stokes' theorem $\int_{\sigma}d\omega = \int_{\partial\sigma}\omega$. The nilpotent property $D^I_{k+1}D^I_k = 0$ follows directly from the fact that the boundary of a boundary always vanishes, i.e., $\partial\partial = 0$.

The discrete Hodge star $S^I_k$ is a diagonal matrix with each diagonal entry given by the ratio between the volume of the dual $(m\!-\!k)$-cell and a primal $k$-cell $\sigma_k$, that is, $l^{m-k}/l^k = l^{m-2k}$. It is induced by the Hodge star operator $\star$ through local averaging
\begin{align}
	\frac{1}{\lvert\sigma_k\rvert}\int_{\sigma_k}\omega \approx \frac{1}{\lvert\star\sigma_k\rvert}\int_{\star\sigma_k}\star\omega,
\end{align}
which establishes a one-to-one correspondence between the primal $k$-forms on $I_m$ and their dual $(m\!-\!k)$-forms on the dual grid. The discrete codifferential is assembled from the discrete differential and Hodge star operators as $\delta^I_k = (S^I_{k-1})^{-1}(D^I_{k-1})^TS^I_k$, following its smooth counterpart.

To construct the discrete weighted differential operators, we let $f$ be a function defined on the entire grid $I_m$ and discretize it using the following strategy: for each degree $k$, we let $f^I_k$ be the diagonal matrix with each diagonal entry computed as the average of $f$ over the vertices of the corresponding primal $k$-cell. Denote by $W^I_k = \exp(f^I_k)$ its exponential. The discrete weighted Hodge star is then given as $S^I_{k, f} = (W^I_k)^{-1}S^I_k = S^I_k(W^I_k)^{-1}$, where $S_k$ and $W^I_k$ are both diagonal and therefore commute. By computation, we obtain a discretization of the weighted codifferential $\delta_f = e^f\delta e^{-f}$ as follows:
\begin{align}
    \delta^I_{k, f} = (S^I_{k-1, f})^{-1}(D^I_{k-1})^TS^I_{k, f}.
\end{align}
The corresponding weighted $L^2$ inner product for two discrete differential forms is then given by:
\begin{align}\label{eq:WeightedL2InnerProduct.discrete}
    (A_k, B_k)^I_f = A_k^TS^I_{k,f}B_k.
\end{align}
Note that discretizing the weighted Hodge Laplacian $\Delta^k_f$ leads to a nonsymmetric matrix. Analogous to the unweighted case \cite{zhao20193d}, we instead consider the discrete counterpart of $\star_f\Delta_f$, which is of the form
\begin{align}
    L^I_{k, f} = (D^I_k)^TS^I_{k+1, f}D^I_k + S^I_{k, f}D^I_{k-1}(S^I_{k-1, f})^{-1}(D^I_{k-1})^TS^I_{k, f}.
\end{align}
Here all operators are taken to be null for $k<0$ and $k>m$. The discrete weighted Hodge Laplacian $L_{k, f}$ is positive semidefinite symmetric, and thus its eigenvalues are non-negative.

\subsection{Discretization on $M$}

We now construct discrete differential operators on $M$, which is represented as a sublevel set of a level set function on a regular Cartesian grid. Note that the Hodge duality identifies the two types of boundary conditions. Here we present only the discretization of the weighted framework under normal boundary conditions. The tangential case follows similarly. Unless stated otherwise, all operators are understood to be under normal boundary conditions.


To restrict all computations to relevant cells and impose normal boundary conditions, we use projection matrices through the inclusion or exclusion of the entire $k$-cells. We consider the normal support, which is defined to be the set of all primal cells with at least one vertex inside $M$, following \cite{su2024persistent}. Denote by $P_{k}$ the projection matrix that maps all $k$-cells of the entire grid $I_m$ to the normal support. The matrix $P_k$ can be obtained from the identity matrix by eliminating the rows corresponding to $k$-cells outside the support.

The discretization of the Hodge star operators on $M$ also requires incorporating the normal boundary conditions. For each diagonal entry in the discrete Hodge star matrix, we replace the full volume of the primal $k$-cell by the $k$-volume of its intersection with $M$ and keep the dual cell volumes unchanged. In addition, we apply a small numerical perturbation to the level set function evaluated at primal/dual grid points to have an absolute value above $\epsilon = 10^{-5}l$, which ensures the fractional $k$-volumes are well-behaved. Denote by $\tilde S^I_{k}$ the resulting diagonal Hodge star matrix on the entire grid $I_m$. Then we have the following set of differential operators defined on the normal support of $M$:
\begin{align}
D_{k} = P_{k+1} D^I_k P_{k}^T, \qquad
S_{k} = P_{k} \tilde S^I_{k} P_{k}^T.
\end{align}
The nilpotent property $D_{k+1}D_{k} = 0$ remains valid \cite{su2024persistent}.

To incorporate weights in the discrete operators, we let $f$ be a function on $M$. Denote by $f_{k}$ its discretization on the normal support for degree $k$, which is diagonal. We compute the diagonal entry of $f_{k}$ simply by taking the average of $f$ evaluated at the associated grid points of the corresponding primal $k$-cell inside $M$, along with the intersection points if the cell intersects the boundary. Let $W_{k} = \exp(f_{k})$. The discrete weighted Hodge star on $M$ is then $S_{k, f} = W_k^{-1}S_k = S_kW_k^{-1}$ and the discrete weighted codifferential is $\delta_{k, f} = (S_{k-1, f})^{-1}(D_{k-1})^TS_{k, f}.$ Note that we omitted the potential sign change $(-1)^k$ in the continuous counterpart as it does not affect the Laplacian construction. In addition, we have the weighted $L^2$-inner product for two discrete differential forms on the normal support as follows 
\begin{align}
    (A_k, B_k)_f = A_k^TS_{k, f}B^k.
\end{align}
The discrete weighted Hodge Laplacian on $M$ can be assembled as
\begin{align}
    L_{k, f} = D_k^TS_{k+1, f}D_k + S_{k, f}D_{k-1}(S_{k-1, f})^{-1}(D_{k-1})^TS_{k, f}.
\end{align}
Consistent with the smooth theory, the kernel of the discrete weighted Hodge Laplacian $L_{k, f}$ is fully determined by the manifold topology, and its kernel dimension is given by the Betti number $\beta_{m-k}$, independent of $f$. Note that the statement is for normal boundary conditions. The kernel of the discrete tangential weighted Hodge Laplacian (not discussed here) in degree k has dimension $\beta_{k}$.

The spectra of $S_{k,f}^{-1}L_{k,f}$ carry rich information about the topological and geometric information about the underlying manifold. Besides the null space that captures the manifold topology, the nonzero part reflects geometry. It is known that the first non-zero eigenvalue, called the Fiedler value, describes connectivity, and eigenvalue multiplicities often reflect symmetries of the shape of data. In addition, the weighted setting allows the study of features across different weights. In particular, one is able to highlight and differentiate the local geometry by selecting weights $f$ that emphasize regions of interest. 

\begin{rem}
In the case that $f$ is constant, the matrix $f_k$ is proportional to the identity matrix with the same scaling factor across all degrees, and so is the matrix $W_k$. The discrete weighted codifferential $\delta_{k,f}$ then reduces to the standard codifferential $\delta_k$, and the discrete weighted Hodge star $S_{k,f}$ is identical to the standard discrete Hodge $S_k$ up to a constant. It follows that the discrete weighted $L^2$-inner product~\eqref{eq:WeightedL2InnerProduct.discrete} coincides with the discrete Hodge $L^2$-inner product. In addition, the discrete weighted Laplacian $L_{k, f}$ is identical to the standard discrete Hodge Laplacian, obtained by replacing $S_{k,f}$ with $S_{k}$. Thus, $S_{k,f}^{-1}L_{k,f}=S_k^{-1}L_k,$ and
the discrete weighted framework recovers the unweighted case when $f$ is constant, analogous to the smooth setting.
\end{rem}

\begin{rem}
The discrete weighted Hodge Laplacians on normal support and tangential support play a crucial role in the implementation of the Hodge decomposition~\eqref{eq.hd.5subspaces} for compact domains on Cartesian grids in 2D and 3D. They can be used to compute the potentials associated with the decomposed components. Furthermore, the associated eigenvectors of $0$ eigenvalues of the Laplacian for each $k$ on a chosen support form a basis of the space of discrete harmonic fields with the corresponding boundary conditions.
\end{rem}

Note that the discrete Hodge star defined on a Cartesian grid is almost identical to a rescaled identity matrix. One may therefore replace the Hodge star with the identity matrix to simplify the computations. This leads to the BIG Laplacians \cite{ribando2024combinatorial}, which were introduced to facilitate the comparison and contrast between Hodge Laplacians and combinatorial Laplacians. The BIG Laplacians have the same rank deficiencies as the corresponding discrete Hodge Laplacians, and thus capture the correct cohomology. In addition, their spectra converge to those of Hodge Laplacians up to a scaling value. 

Analogous to the unweighted case, we define the weighted BIG Laplacian by replacing the Hodge star $S_k$ with the identity matrix. The weighted Hodge star is then given as the inverse of the exponential weight matrix $S_{k, f} = (W_k)^{-1}$. This finally leads to the weighted BIG Laplacian as follows:
\begin{align}\label{eq:BIGLaplacian.weighted}
    L^B_{k, f} = D_k^TW_{k+1}^{-1}D_k + (W_k)^{-1}D_{k-1}W_{k-1}D_{k-1}^{T}(W_k)^{-1}.
\end{align}
In the case that $f$ is constant, the weighted BIG Laplacian reduces to the BIG Laplacian up to a constant. 
Note that the rank deficiency of the discrete Laplacian depends on the differential and projection matrices. The kernel size of the weighted BIG Laplacian is unchanged as in the unweighted setting, given by the same Betti number. In practice, the weighted BIG Laplacians can be used in the same way as the weighted Hodge Laplacians to study the topological and geometric information of the underlying manifold.

\subsection{Discrete spectral analysis of Laplacians}

The eigenvalue analysis of discrete Laplacians in the unweighted setting has been discussed in detail in \cite{su2024persistent}. The situation for the weighted framework remains the same, as incorporating the weight function modifies only the diagonal entries of the discrete Hodge stars and does not change the boundary conditions. Therefore, we omit the details and refer the reader to \cite{su2024persistent} for the unweighted case.

The eigenvalues and eigenvectors of a discrete Laplacian $L_{k, f}$ can be solved by considering the generalized eigenvalue problem
\begin{align}\label{eq:generalized_eig_problem}
    L_{k,f}V = \lambda S_{k, f}V,
\end{align}
where $L_{k,f}$ can be interpreted as either the weighted Hodge Laplacian or the weighted BIG Laplacian (with $S_{k, f}$ set to identity) under one type of the boundary conditions, $\lambda$ is an eigenvalue and $V$ is the associated eigenvector, both depending on $f$. Let $\bar L_{k, f} = S_{k, f}^{-1/2}L_{k, f} S_{k, f}^{-1/2}$, and $\bar V = S_{k, f}^{1/2}V$. The formulas~\eqref{eq:generalized_eig_problem} can then be written as a regular eigenvalue problem:
\begin{align}
\bar L_{k, f} \bar V = \lambda \bar V.
\end{align}
For our implementation, $S_{k,f}$ is diagonal, and thus its square root and inverse square root are also diagonal and efficient to evaluate. So the spectral analysis of the modified Laplacian $\bar L_{k, f}$ can be efficiently converted to the study of the singular spectrum of the discrete operator $\bar D_{k, f} = S_{k+1, f}^{1/2}D_kS_{k, f}^{-1/2}$. Specifically, the spectrum of $\bar L_{k, f}$ on the normal support is given by the union of the squared nonzero singular values of $\bar D_{k-1, f}$ and $\bar D_{k, f}$, together with $0$ of multiplicity $\beta_{m-k}$, i.e., the $(m\!-\!k)$-th Betti number.

In the case that $\dim(M) = 3$, there are eight discrete weighted Hodge Laplacians for $k=0,1,2,3$ under the normal and tangential boundary conditions. Due to the duality of these two types of boundary conditions, the study of the spectra of these Laplacians, as noted above, can finally be reduced to the study of the singular values of the three discrete differentials $\bar D_{0, f}$, $\bar D_{1, f}$ and $\bar D_{2, f}$ with one type of the boundary conditions.


\section{Protein Flexibility Analysis}
\label{sec:proteinFlexibilityAnalysis}
To validate the effectiveness of our framework for extracting the topological and geometric features of data, and in particular for highlighting and differentiating local features, we apply it to protein flexibility analysis, also called B-factor analysis. The protein flexibility characterizes the dynamic behavior of atoms or amino acids within a protein structure, and is quantified by the B-factor (also known as the Debye-Waller factor or temperature factor), which measures an atom's displacement from its mean position. It provides valuable information about protein's structural and thermal stability, activity of the region, and other protein functions, and has important applications such as docking and structure-based drug design \cite{chandrika2009managing,carlson2000accommodating,teague2003implications}.

Many frameworks have been developed for protein flexibility analysis, including molecular dynamics (MD) \cite{mccammon1977dynamics},
normal-mode analysis (NMA) \cite{brooks1983charmm,go1983dynamics,levitt1985protein,tasumi1982normal},
elastic network models (ENMs) \cite{tirion1996large} with their variants Gaussian network model (GNM) \cite{bahar1998vibrational,bahar1997direct} and
anisotropic network model (ANM) \cite{atilgan2001anisotropy}, the flexibility and rigidity index (FRI) methods \cite{xia2013multiscale, opron2014fast} and more recently, blind prediction machine learning models that utilize topological and geometric features of data, such as multiscale differential geometry learning (mDGL) model \cite{feng2025multiscale}, persistent sheaf Laplacian (PSL) model \cite{hayes2025persistent} and commutative algebra learning (CAL) model \cite{zhang2026commutative}. These machine learning models mDG, PSL and CAL have been reported to achieve higher accuracy than previous models.









To demonstrate the performance of our framework, we follow the same evaluation protocols used by the mDG, PSL and CAL models, including the same 10-fold cross-validation and data splitting strategy on the same dataset, and compare our results with theirs. We built two models, one using only features extracted from the weighted Hodge Laplacians, and a consensus model that incorporates additional features from PDB files and those generated using STRIDE, a software that generates protein secondary structure assignment from atomic coordinates. The dataset is from Opron et al.~\cite{opron2014fast, opron2015communication}, which contains a total of 364 proteins. In the experiments, several proteins are excluded due to nonphysical issues (0 B-factor values) and inconsistencies with STRIDE, see \cite{feng2025multiscale} for the list of excluded PDB IDs. Finally, we consider only a total of $346$ proteins. Data is available at \url{https://github.com/fenghon1/MDG_bfactor}.

Utilizing the 3D coordinates of atoms in proteins, we generate a manifold with boundary for each protein based on a discrete-to-continuum mapping. We then compute for each atom in the protein the weighted Hodge Laplacian with an atom-specific weight function defined on the resulting manifold. The eigenvalues of the weighted Hodge Laplacian, termed WHL features, are then used as a descriptor of the atom, which includes the global topological feature of the manifold coming from the number of 0 eigenvalues, i.e., the 0th Betti number $\beta_0$, and also the atom-specific geometric features from the non-zero eigenvalues. The associated B-factor serves as the label for an atom. The 3d coordinates and B-factors of atoms are both recorded in PDB files. As analyzing B-factors of $C_{\alpha}$ atoms is standard for describing backbone mobility and dynamics resulting from protein dynamics, we consider only $C_{\alpha}$ atoms in the experiments. The feature extraction procedure of our proposed weighted Hodge Laplacian learning (WHLL) framework is illustrated in Fig.~\ref{fig.mfld.1IDP}.

In the following, we present details of the feature extraction algorithm for the protein data, including the generation of manifolds, the construction of atom-specific weight functions, and details of our blind prediction machine learning algorithms. Finally, we report our results to demonstrate the effectiveness of our models. The computation of eigenvalues of Laplacians is implemented in MATLAB, while the machine learning models are employed using a gradient boosting regressor (GBR) module from Scikit-learn.

\subsection{Manifold generation}
\label{sec:manifold_generation}


\begin{figure}[t]
	\centering
	\includegraphics[width=.95\linewidth]{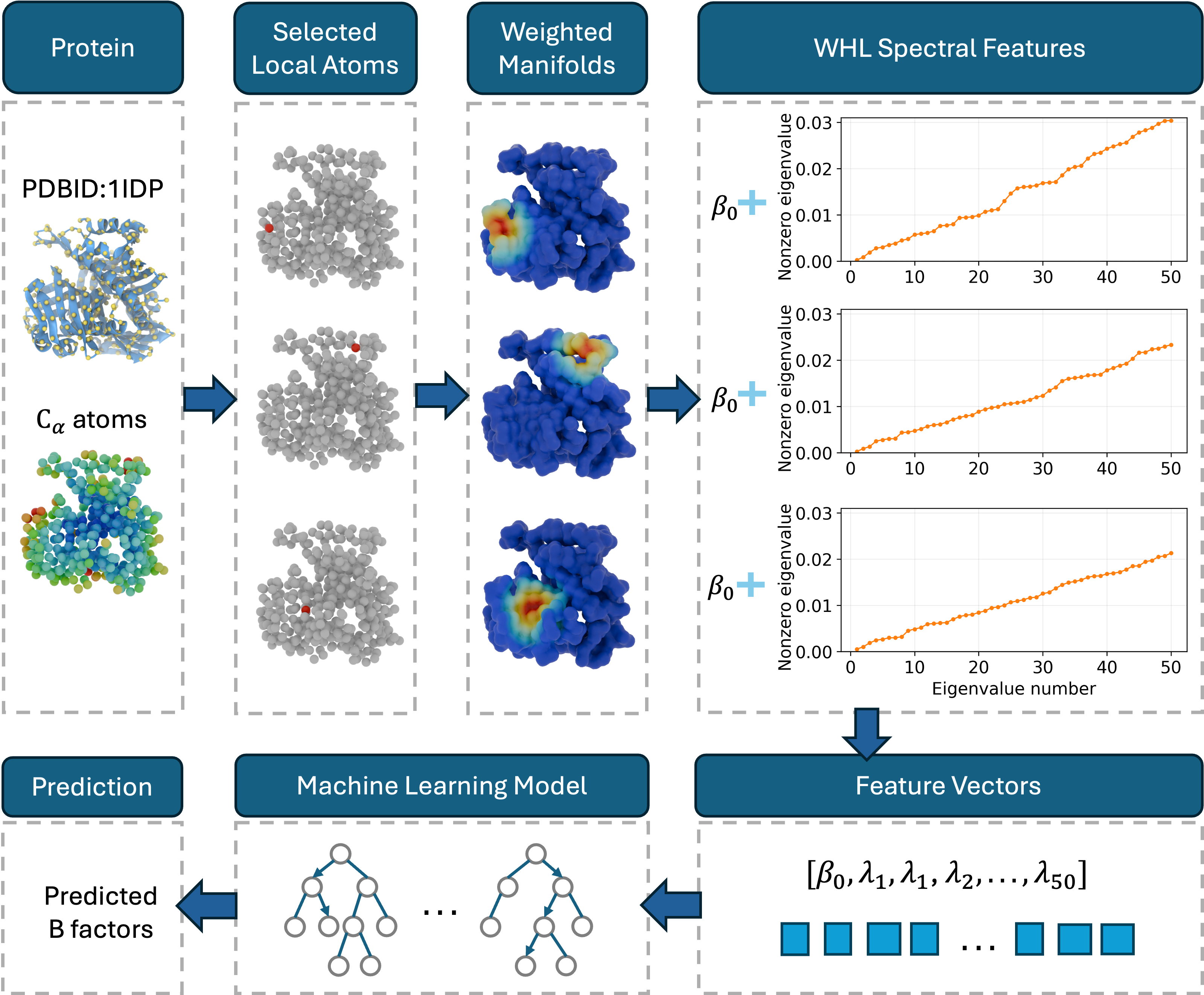}
	\caption{Schematic illustration of the weighted Hodge Laplacian learning (WHLL) framework for atom-wise B-factor prediction. First row, first column: Protein 1IDP and its $C_{\alpha}$-atoms colored by a rainbow colormap (blue\texttt{->}cyan\texttt{->}green\texttt{->}yellow\texttt{->}orange\texttt{->}red), where blue denotes the lowest B-factors and red the highest. First row, second column: Three selected $C_{\alpha}$-atoms. First row, third column: Manifolds generated from \eqref{eq.lvf.rho} with atom-specific weights~\eqref{eq:weight} that emphasize local regions of the $C_\alpha$ atoms shown in the second column. Red indicates higher weights, and dark blue represents near-zero weights. First row, fourth column: The resulting features for these atoms consisting of the 0-th Betti number $\beta_0$ and the first $k$ nonzero eigenvalues, which are used as inputs for the machine learning models (with $k=50$ for illustration). Second row: Machine learning model for predicting atom-wise B-factors.}
	\label{fig.mfld.1IDP}
\end{figure}

Note that the original data are point clouds given as 3D coordinates of atoms for each protein. We need to generate a manifold representation for each protein so that the data are compatible with our weighted Hodge Laplacian framework. This process can be realized by a discrete-to-continuum mapping using the flexibility and rigidity index \cite{nguyen2019dg}. 

Denote by $\{\mathbf{x_1}, \mathbf{x_2}\ldots, \mathbf{x_s}\}$ the 3D coordinates of all $C_\alpha$ atoms in a protein. We define a discrete-to-continuum mapping as follows, given by the negative sum of Gaussian densities at all positions of these atoms:
\begin{align}\label{eq.lvf.rho}
	\rho(\mathbf{x}, \tau) = -\sum^{s}_{i = 1}\exp\left(-\left(\frac{\norm{\mathbf{x} - \mathbf{x_i}}}{\tau r}\right)^2\right),
\end{align}
where $\norm{\mathbf{x} - \mathbf{x_i}}$ is the Euclidean distance from $\mathbf{x}$ to the position $\mathbf{x_i}$ of the $i$-th atom, $\tau$ is a scalar value, and $r=1.7$ is the van der Waals radius of the $C_\alpha$ atom. Given an isovalue $c$, the sublevel set
\begin{align}
	M = \{\mathbf{x}\,\vert\, \rho(\mathbf{x}, \tau)\leq c \}
\end{align}
then leads to a compact manifold in $\RR^3$ with its boundary given by the isosurface $\partial M = \{\mathbf{x} \, \vert\, \rho(\mathbf{x}, \tau) = c \}$ for the protein. In our experiments, we set $\tau = 1, c = -0.1$ to generate manifolds for all proteins. These values are chosen to produce smooth manifolds with stable topology and to avoid numerical issues. In Fig.~\ref{fig.mfld.1IDP}, we present an example of the resulting manifolds for three selected $C_{\alpha}$ atoms in protein 1IDP using \eqref{eq.lvf.rho}.

Note that the level set function~\eqref{eq.lvf.rho} is a special case of the flexibility and rigidity index (FRI) density function \cite{nguyen2019dg}, which is known to be stable for converting point clouds into continuous embeddings and has been used to generate protein boundary surfaces and interactive manifolds \cite{chen2021evolutionary,nguyen2019dg}. Other choices of FRI density functions can also be employed to generate the manifold for each protein.




\subsection{Atom-specific weight construction}
\label{sec:weight_construction}

Notice that the B-factor of an atom is influenced by the local packing and interactions of nearby atoms. To characterize and differentiate its local feature with position $\mathbf{\bar x}$ in a protein, we consider only atoms within its neighborhood with a cutoff distance $d = 11$ from $\mathbf{\bar x}$, and use them to construct the atom-specific weight function for the atom. This choice of cutoff distance has been confirmed by our experiments that the local atomic structures can be captured in the weight functions, and the resulting spectra of Laplacians distinguish between atoms. 

Given an atom in a protein with position $\mathbf{\bar x}$, let $N$ be the number of atoms in the neighborhood within the cutoff distance $d$ of $\mathbf{\bar x}$. We define a weight function for this atom by the following distance-weighted sum of densities 
\begin{align}\label{eq:weight}
    f(\mathbf{x}) = \sum^{N}_{j = 1}\omega_j\exp\left(-\left(\frac{\norm{\mathbf{x} - \mathbf{x_j}}}{\tau r}\right)^2\right),\quad \omega_j = \exp\left(\frac{-\norm{\mathbf{x_j} - \mathbf{\bar x}}^2}{2\eta^2}\right),
\end{align}
where $r=1.7$ is the van der Waals radius of the $C_\alpha$ atom, $\tau = 5$ and $\eta=4$ are chosen in our experiments to include wider local interactions of the atom and to avoid overemphasizing its very near neighbors. Examples of weight functions that emphasize local regions of three different $C_{\alpha}$ atoms in protein 1IDP are shown in Fig.~\ref{fig.mfld.1IDP}.

The formulation of the weight function \eqref{eq:weight} ensures a decay of influence with distance, which is consistent with the nature of atomic interactions in a molecular system. The parameters $\tau$ and $\eta$ control the locality and the decay rate, allowing our model to emphasize atom-specific local structures while limiting the influence of distant atoms.

For all $C_{\alpha}$ atoms in each protein, we compute the atom-specific weight~\eqref{eq:weight}, and then their weighted Hodge Laplacians on the manifold of the protein, generated following Sec.~\ref{sec:manifold_generation}. Their kernels are of the same dimension, independent of the choice of $f$ and determined by the manifold topology. Their non-zero spectra vary with weights for different atoms, highlighting and differentiating their local geometric features.

\subsection{Feature extraction}
\label{sec:feature_extraction}

With the manifold representation and the atom-specific weight functions at hand, we implement our framework on a regular Cartesian grid with the chosen normal boundary conditions. For all atoms across proteins, we fix a Cartesian grid with grid length $1$. This pre-chosen Cartesian grid ensures that all computation is consistent and the resulting features are comparable for all atoms regardless of the protein IDs, as they are obtained following the same discretization scheme. The choice of grid length $1$ provides sufficient grid resolution for accurate computation of Laplacians so that no topological information is lost due to numerical errors caused by low resolution. 
All discrete operators are assembled following Sec.~\ref{sec:discretization}.

In our experiments, we compute only the weighted BIG Laplacian $L^B_{3, f}$ under normal boundary conditions for each atom due to the computational efficiency. The spectra of $L^B_{3, f}$ serve the same role as the $0$-th Laplacian in the combinatorial setting and have proven effective and successful in many machine learning tasks \cite{cang2018integration, liu2017forging, meng2021persistent, wang2020persistent, su2024persistent}. The number of the $0$ eigenvalues of the Laplacian provides the $0$-th Betti number $\beta_0$, i.e., the number of connected components of the manifold, while the non-zero eigenvalues encode the local geometric structures induced by the atom-specific weights. We take the $0$-th Betti number $\beta_0$, together with the first $k$ non-zero eigenvalues. This leads to a total of $k\!+\!1$ features for each atom, which are then fed into our machine learning models. As demonstrated in our experiments, these features are sufficient to validate our frameworks in capturing the topological and geometric features and differentiating atoms.

In addition to the WHL features, we incorporate $3$ global and $9$ local protein features to build consensus machine learning models. The $3$ global features are the R-value, protein resolution, and the number of heavy atoms extracted from the PDB files. The $9$ local features are obtained using STRIDE \cite{heinig2004stride}. These features include three packing density values with cutoff distances (short $d\!<\!3$, medium $3\!\leq\! d<\!5$, and long-ranged $d\!\geq\!5$), amino acid type, occupancy, and secondary structure type, $\phi$ and $\psi$ angles, and residue solvent-accessible area. For the consensus model, each atom has $(k\!+\!1)\!+\!12$ features with $k\!+\!1$ being the WHL features and $12$ the additional global and local features.

\subsection{Machine learning algorithm}


To evaluate the performance of our models, we carry out the 10-fold cross-validation under two setups: protein-level and atom-level. At the protein level, we split in each fold the set of 346 proteins into training and test subsets so that all $C_{\alpha}$ atoms from a given protein are in the same subset. At the atom level, we split the set of all $C_{\alpha}$ atoms (more than 74,000 atoms) across proteins into training and test subsets. For each setup, the process is repeated over all 10 folds, and we report the average performance.

For both the WHL and consensus models, we employ the GBR module from Scikit-learn 1.4.2 for the blind prediction of B-factors of atoms, with the following GBR parameters: \texttt{n\_estimators} = \texttt{1{,}000}, \texttt{max\_depth} = \texttt{7}, \texttt{min\_samples\_split} = \texttt{5}, \texttt{learning\_rate} = \texttt{0.002}, \texttt{loss} = \texttt{squared\_error}, \texttt{subsample} = \texttt{0.8}, and \texttt{max\_features} = \texttt{sqrt}. The Pearson correlation coefficient (PCC) is used as the evaluation metric in our models.

\subsection{Experimental results}


\begin{table}[htbp]
\centering
\caption{Average Pearson correlation coefficients (PCC) for the WHL model and the consensus model, computed via 10-fold cross-validation at the protein and atom levels on 346 proteins with more than 74,000 $C_{\alpha}$ atoms using the gradient boosting regressor (GBR) module in scikit-learn. Results are compared with mDGL, PSL and CAL.}
\label{table:PCCresults}

\begin{threeparttable}		
\begin{tabularx}{0.8\linewidth}{l *{5}{S}}
	\toprule
	& {mDGL} & {PSL} & {CAL} & {WHL} & {WHL\_consensus} \\
	\midrule
	Protein-level & 0.407 & 0.452 & 0.456 & 0.480 & {\bf 0.524} \\
	Atom-level    & 0.859 & 0.840 & 0.855 & 0.842 & {\bf 0.862} \\
	\bottomrule
\end{tabularx}
\vspace{0.25em}
\footnotesize
\begin{tabular}{@{}r@{:\ }l@{}}
	mDGL & Multiscale Differential Geometry Learning \cite{feng2025multiscale}\\
	PSL  & Persistent Sheaf Laplacians \cite{hayes2025persistent}\\
    CAL  & Commutative Algebra Learning \cite{zhang2026commutative}\\
	WHL  & Weighted Hodge Laplacians
\end{tabular}
\end{threeparttable}
\end{table}

As presented in Sec.~\ref{sec:feature_extraction}, the number of WHL features for each atom is given by $k\!+\!1$, consisting of the $0$-th Betti number $\beta_0$ and the first $k$ non-zero eigenvalues of the weighted Hodge Laplacian defined on the corresponding manifold. To find the optimal $k$ that leads to the best performance of our predictive modules, we vary $k$ with an increment of $10$, i.e, $k = 10, 20, \ldots$ and perform the 10-fold cross-validation for each value of $k$. The results suggest that at the protein level, the optimal PCC value can be achieved by the WHL model when $k=120$ (PCC = 0.480), and by the consensus model when $k=90$ (PCC = 0.524). At the atom level, the WHL model is optimal when $k = 150$ (PCC = 0.842), while the consensus model is optimal when $k=70$ (PCC = 0.862). In Table~\ref{table:PCCresults}, we report these average PCC values for the WHL and consensus models at both the protein level and the atom level, along with the performance of mDG, PSL and CAL. The best performance is achieved when using the consensus model, leading to a PCC value of $0.524$ at the protein-level and a PCC value of $0.862$ at the atom level, both higher than previous models. The results demonstrate the effectiveness and also the promise of our models in capturing the topological and geometric features of data, and their capability to highlight and differentiate local features.

\section{Conclusions}
\label{sec:conclusion}

In this work, we develop a weighted de Rham-Hodge Laplacian framework for manifolds with boundary. The theory extends the classical setting by incorporating a smooth function on the manifold and is obtained by modifying only the codifferential operator. Under appropriate boundary conditions, the kernels of the resulting weighted Hodge Laplacians (or drifting Hodge Laplacians) coincide with the weighted harmonic spaces, and are isomorphic to the corresponding manifold cohomology groups, as in the unweighted case. 
Their harmonic spectra therefore capture the global topological information about the underlying manifolds, while the non-harmonic spectra encode geometric structure with respect to a chosen weight. In addition, this method allows one to highlight and differentiate the local geometric features that are inaccessible in the unweighted case. Numerically, we provide a DEC-based computational framework for the weighted Hodge Laplacians on regular Cartesian grids and also introduce the weighted boundary-induced graph (BIG) Laplacians without Hodge star to facilitate computations. We validate the effectiveness of our framework by building blind prediction machine learning models for protein flexibility analysis, which requires localized information. Results show the promise of our framework for analyzing the local properties of data on manifolds.

Several directions remain open for further investigation. For example, it is natural to establish a $5$-component Hodge decomposition in the weighted setting for manifolds with boundary, and to develop the corresponding discretization framework. In addition, one can also modify both the differential and codifferential and use the resulting Witten Laplacian for extracting the topological and geometric features. In this context, all standard results, including the adjointness of the modified differential and codifferential, the self-adjointness of the Witten Laplacian, and the corresponding Hodge decompositions, are preserved with respect to the standard Hodge $L^2$-inner product. Therefore, one no longer needs to redefine the $L^2$-inner product, which makes the framework theoretically simple. By introducing weights, these methods are well-suited for data-driven analysis of data with local features. Finally, the proposed method has potential applications to manifold topological learning \cite{su2024persistent}, partial differential equation (PDE)-based learning models, and generative AI with weighted Hodge diffusion. 

\section*{Acknowledgements}
The work of Wei was supported in part by NIH grant R35GM148196, the University of Georgia,  and Georgia Research Alliance. 
The authors thank Professor Shmuel Weinberger for useful discussions.  



\bibliography{refs}

\end{document}